\documentclass{elsart}

\usepackage{natbib}

\usepackage{amssymb,amsmath,  amssymb,mathrsfs}

\numberwithin{equation}{section}
\newtheorem{theorem}{Theorem}
\newtheorem{remark}{Remark}
\newtheorem{lemma}{Lemma}
\newtheorem{corollary}{Corollary}
\numberwithin{theorem}{section} \numberwithin{remark}{section}
\numberwithin{lemma}{section} \numberwithin{corollary}{section}

\journal{Journal of Statistical Planning and Inference}
\begin{document}

\begin{frontmatter}



\title{Locally Most Powerful Sequential Tests of a Simple Hypothesis vs One-Sided Alternatives}


\author[an]{Andrey Novikov}, \ead{an@xanum.uam.mx}
\ead[url]{http://mat.izt.uam.mx/profs/anovikov/en} \author[pn]{Petr Novikov}\ead{pnovi@mail.ru}
\address[an]{Department of Mathematics,\\ Autonomous Metropolitan University - Iztapalapa,\\
San Rafael Atlixco 186, col. Vicentina\\
C.P. 09340, Mexico City, Mexico}
\address[pn] {Department of Mathematical Statistics, Kazan State University,\\
Kremlevskaya 18, Kazan, Russia}
\begin{abstract}
Let $X_1,X_2,\dots$ be a discrete-time stochastic
process with a distribution $P_\theta$, $\theta\in\Theta$, where $\Theta$ is an open subset of the real line. We consider the problem of testing
a simple hypothesis $H_0:$ $\theta=\theta_0$ versus a composite alternative
$H_1:$ $\theta>\theta_0$, where $\theta_0\in\Theta$ is some fixed point. The main goal of this article is to characterize the
structure of locally most powerful sequential tests in this problem.

For any sequential test $(\psi,\phi)$ with a (randomized) stopping rule $\psi$ and a (randomized) decision rule $\phi$ let $\alpha(\psi,\phi)$ be the type I error probability, $\dot \beta_0(\psi,\phi)$ the derivative, at $\theta=\theta_0$, of the power function, and $\mathscr N(\psi)$ an average sample number  of the test $(\psi,\phi)$. Then we are concerned with the problem of maximizing $\dot \beta_0(\psi,\phi)$ in the class of all sequential tests such that
$$
\alpha(\psi,\phi)\leq \alpha\quad\mbox{and}\quad \mathscr N(\psi)\leq \mathscr N,
$$
where $\alpha\in[0,1]$ and $\mathscr N\geq 1$ are some restrictions. It is supposed that $\mathscr N(\psi)$ is calculated under some fixed (not necessarily coinciding with one of $P_\theta$) distribution of the process $X_1,X_2\dots$.

The structure of optimal sequential tests is characterized.
\end{abstract}


\end{frontmatter}

\section{Introduction}
\label{intro}
Let $X_1,X_2,\dots, X_n, \dots$ be a discrete-time stochastic
process with a distribution $P_\theta$, $\theta\in\Theta$, where $\Theta$ is an open subset of the real line. We consider the problem of testing
a simple hypothesis $H_0:$ $\theta=\theta_0$ versus a composite alternative
$H_1:$ $\theta>\theta_0$, where $\theta_0\in\Theta$ is some fixed point. The main goal of this article is to characterize the
structure of locally most powerful, in the sense of \cite{berk}, sequential tests in this problem.

We follow  \cite{novikov09b} in the definitions and notation related to sequential hypothesis tests, as well as their interpretation and characteristics (see also \cite{wald50}, \cite{ferguson}, \cite{degroot}, \cite{schmitz}, \cite{ghosh}, among many others).

In particular, we say that a pair $(\psi,\phi)$ is a sequential hypothesis test
if
$$\psi=\left(\psi_1,\psi_2, \dots ,\psi_{n},\dots\right)\quad\mbox{and}\quad \phi=\left(\phi_1,\phi_2, \dots ,\phi_n,\dots\right),$$
where  the  functions $$\psi_n=\psi_n(x_1, x_2,\dots, x_n)
\quad\mbox{and}\quad \phi_n=\phi_n(x_1, x_2,\dots, x_n)$$
 are supposed to be  measurable functions with values in $[0,1]$, $n=1,2,\dots$.

For any stage $n=1,2,\dots$, the
value of $\psi_n(x_1,\dots,x_n)$ is interpreted as the conditional
probability {\em to stop and proceed to decision making}, given that
the experiment came to stage $n$ and that the observations of
the process up to  this stage were $(x_1, x_2, \dots,
x_n).$ The rules $\psi_1,\psi_2,\dots$ are successively  applied until the experiment eventually stops.

It is supposed that when the experiment stops, at some stage $n\geq 1$, the decision rule
$\phi_n$ will be applied to make a decision.
The value of
$\phi_n(x_1,\dots, x_n)$ is interpreted  as the conditional
probability {\em to reject} the null-hypothesis $H_0$, given that
the data observed up to this stage, were $(x_1,\dots, x_n)$.

The stopping rule $\psi$ generates, by the above process, a random
variable $\tau_\psi$ ({\em stopping time}) whose distribution is
given by
\begin{equation}\label{1.1}
P_\theta(\tau_\psi=n)=E_\theta(1-\psi_1)(1-\psi_2)\dots (1-\psi_{n-1})\psi_n.
\end{equation}

Here, and throughout the paper, $E_\theta(\cdot)$ stands for the expectation with respect to the distribution $P_\theta$ of the process $X_1, X_2, \dots$.

In (\ref{1.1}), we suppose that $\psi_n=\psi_n(X_1,X_2,\dots,X_n)$, unlike its previous definition as $\psi_n=\psi_n(x_1,x_2,\dots, x_n)$.
We do this intentionally and systematically
throughout the paper, applying, generally, for any $F_n=F_n(x_1,x_2,\dots,x_n)$ or $F_n=F_n(X_1,X_2,\dots,X_n)$, the following rule:
if $F_n$ is under the probability or
expectation sign, then it is $F_n(X_1,\dots, X_n)$, otherwise it is
$F_n(x_1,\dots, x_n)$.

To characterize the duration of the sequential experiment, the {\em average sample number} is used:
\begin{equation}\label{1.2}
{\mathscr N}(\psi)=E\tau_\psi=\begin{cases}\sum_{n=1}^\infty nP(\tau_\psi=n),\; \mbox{if}\; P(\tau_\psi<\infty)=1,\cr
\infty,\; \mbox{ otherwise.} \end{cases}
\end{equation}

For a sequential test $(\psi,\phi)$ let us define its {\em power function}
at $\theta$ as
\begin{equation}\label{1.3}
\beta_\theta(\psi,\phi)=P_{\theta}(\mbox{reject}\; H_0)=\sum_{n=1}^\infty E_\theta(1-\psi_1)
\dots(1-\psi_{n-1})\psi_n\phi_{n}.
\end{equation}
The {\em type I error probability} of the test $(\psi,\phi)$ is defined as
\begin{equation*}
\alpha(\psi,\phi)=\beta_{\theta_0}(\psi,\phi).
\end{equation*}

Our main goal is characterizing tests which  maximize the derivative of the power function at $\theta=\theta_0$,
$
\dot\beta_{\theta_0}(\psi,\phi),
$
among all sequential tests $(\psi,\phi)$ such that
\begin{equation}\label{1.4}
\alpha(\psi,\phi)\leq \alpha,
\end{equation}
and
\begin{equation}\label{1.5}
\mathscr N(\psi)\leq \mathscr N,
\end{equation}
where $\alpha\in[0,1)$ and $\mathscr N\geq 1$ are some restrictions. In case this test exist, it is called the locally most powerful test (see \cite{berk}, \cite{roters}).

There is a natural candidate for the distribution under which $\mathscr N(\psi)=E\tau_\psi$  is calculated in (\ref{1.5}): it is  $P_{\theta_0}$ (see \cite{berk} or  \cite{schmitz}).
 Nevertheless, we pose a more general problem in this article, supposing that $E\tau_\psi$
 is calculated under an arbitrary (but fixed) distribution of the process.
 In particular, it may  be useful  to employ as $P$ a ``mixed'' distribution defined as $$P(\cdot)=\int_\Theta P_{\theta}(\cdot)d\pi(\theta),$$ where $\pi$ is some probability measure (see Section 4.2 in \cite{novikov09b} for a good reason for doing so).

\section{Assumptions and Notation}
We suppose throughout the paper that, under $P_\theta$, for all $\theta\in\Theta$,
the vector $(X_1, X_2, \dots , X_n)$ has a
probability ``density" function
\begin{equation*}
f_\theta^n=f_\theta^{n}(x_1, x_2, \dots, x_n)
\end{equation*}
(Radon-Nikodym derivative of its distribution) with respect to a
product-measure $$\begin{array}{ccc}\mu^n&=&\underbrace{\mu\otimes \mu\otimes \dots\otimes\mu},\\
&&n\quad \mbox{times}\end{array}$$
with some $\sigma$-finite measure $\mu$ on the respective space.

We will also suppose that the  distribution $P$ of the process
used for calculating (\ref{1.2}) is some arbitrary (but fixed) distribution
such that $(X_1,\dots,X_n)$ has a ``density''
$
f^n(x_1,\dots,x_n)
$ with respect to $\mu^n$, $n=1,2,\dots$.

The following assumption is basic for the differentiability of power functions, a sort of which is obviously needed in view of the problem formulation in the Introduction.

{\sc Assumption 1.}{\em For any $n\geq 1$ there exists a measurable integrable function $\dot f_0^n$ such that
$$
\int\left|f_{\theta_0+h}^n-f_{\theta_0}^n-h\dot f_0^n\right|d\mu^n=o(h),
$$
as $h\to 0$ holds.}

Assumption 1 is nothing more than the $L_1(\mu^n)$-differentiability  of the joint density function $f_\theta^n$, with respect to $\theta$, at $\theta=\theta_0$, for any $n=1,2,\dots$.

In particular, it follows from Assumption 1 that for any measurable function $\phi_n=\phi_n(x_1,\dots, x_n)$, $0\leq \phi_n\leq 1$,
\begin{equation}\label{2.1}
\int\phi_n \left(f_{\theta_0+h}^n-f_{\theta_0}^n-h\dot f_0^n\right)d\mu^n=o(h),
\end{equation}
as $h\to 0$. In fact, it is easy to see that (\ref{2.1}) is equivalent to Assumption 1.

(\ref{2.1}) means that the power function $\beta_\theta(n,\phi)$ of any fixed sample-size test based on the first $n$ observation, is differentiable at $\theta=\theta_0$, and that its derivative is
$$
\dot \beta_{\theta_0}(n,\phi)=\int \dot f_0^n d\mu^n.
$$
(See Conditions C1 to C3 in \cite{novikov06} and similar conditions, for independent and identically distributed (i.i.d.) observations, in \cite{mueller} in relation to differentiability of power functions.)

It is easy to see that if the partial derivative of $f_\theta^n$ with respect to $\theta$ exists $\mu^m$-almost everywhere at $\theta=\theta_0$, then  it follows from (\ref{2.1}) that \begin{equation}\label{2.2}\dot f_0^n=\left.\frac{\partial f_\theta^n}{\partial \theta} \right|_{\theta=\theta_0},\end{equation}
$\mu^n$-almost everywhere.
In fact, this assumption is used in \cite{berk} along with the condition of differentiability of power function of fixed sample size-tests (see Assumptions 2 and  3 in \cite{berk}), in the i.i.d. case.

The following assumption is needed to treat the optimality in the general case of non-truncated tests below.

{\sc Assumption 2.} {\em The power function of any test $(\psi,\phi)$ such that $E_{\theta_0}\tau_\psi<\infty$, is differentiable at $\theta=\theta_0$, and }
\begin{equation}\label{2.3}
\dot\beta_{\theta_0}(\psi,\phi)=\sum_{n=1}^\infty \int(1-\psi_1)
\dots(1-\psi_{n-1})\psi_n\phi_{n}\dot f_{0}^n\,d\mu^n.
\end{equation}

If the partial derivative (\ref{2.2}) exists, (\ref{2.3}) may be deemed as differentiating across the integral sign, because of (\ref{1.3}).

For i.i.d. observations, there are various conditions which guarantee the differentiability as in Assumption 2 (see, for example, Proposition 1 in \cite{berk}, or related properties in \cite{mueller86} or \cite{irle}).

We will also need the following

{\sc Assumption 3.} {\em There exist $\gamma>0$ and $N_0>0$ such that
\begin{equation}\label{2.4}
E_{\theta_0}\left(\frac{\dot f_0^n}{f_{\theta_0}^n} \right)^2\leq \gamma n
\end{equation}
for all $n\geq N_0$.}

The expectation on the left-hand side of (\ref{2.4}) is the Fisher information contained in $(X_1,\dots, X_n)$. In the i.i.d. case  considered in \cite{berk},  (\ref{2.4}) is obviously an immediate consequence of Assumption 4 \cite{berk}.

To avoid cumbersome notation, we shall  further on  write
$E_0$, $f_0^n$,  $\beta_0$, and $\dot\beta_0$ instead of
$E_{\theta_0}$, $f_{\theta_0}^n$,  $\beta_{\theta_0}$,
and $\dot\beta_{\theta_0}$, respectively.
\section{Reduction to an optimal stopping problem}
\label{reduction}
To proceed with maximizing $\dot\beta_0(\psi,\phi)$ over the tests subject to
(\ref{1.4}) and (\ref{1.5}) let us define the following
Lagrange-multiplier function:
\begin{equation}{\label{reduction5}}
L(\psi,\phi)=L(\psi,\phi;b,c)=c{\mathscr N}(\psi)+b\alpha(\psi,\phi)-\dot\beta_0(\psi,\phi)
\end{equation}
where $c> 0$ and  $b\in \mathbb R$    are some constant
multipliers.

The following theorem is a direct application of the
Lagrange multiplier method to the conditional problem above.

\begin{theorem}\label{t3.1}
Let $\Delta$ be some class  of sequential tests. Let there exist $c> 0$ and $b>0$
 and a test $(\psi,\phi)\in \Delta$ with $L(\psi,\phi;b,c)>-\infty$, such that
\begin{equation}\label{3.2}
L(\psi,\phi;b,c)=
\inf_{(\psi^\prime,\phi^\prime)\in\Delta}L(\psi^\prime,\phi^\prime;b,c)
\end{equation}
and such that \
\begin{equation}\label{3.3}
{\mathscr N}(\psi)={\mathscr N}
\quad\mbox{and}\quad
\alpha(\psi,\phi)=\alpha.
\end{equation}
 Then for any test $(\psi^\prime,\phi^\prime)\in\Delta$ satisfying
 \begin{equation}\label{3.4}
{\mathscr N}(\psi^\prime)\leq{\mathscr N}
 \quad\mbox{and}\quad
 \alpha(\psi^\prime,\phi^\prime)\leq\alpha
 \end{equation}
 it holds
\begin{equation}\label{3.5}
\dot\beta_0(\psi,\phi)\geq  \dot\beta_0(\psi^\prime,\phi^\prime). \end{equation}

 The inequality
in (\ref{3.5}) is strict if at least one of the inequalities in
(\ref{3.4}) is strict.
\end{theorem}

{\em Proof.} It is quite straightforward:

Let $(\psi^\prime,\phi^\prime)\in \Delta$ be any test
satisfying (\ref{3.4}). Because of (\ref{3.3}) and (\ref{3.2}),
\begin{equation*}
   c{\mathscr N}+b\alpha-\dot\beta_0(\psi,\phi)= c{\mathscr N}(\psi)+b\alpha(\psi,\phi)-\dot\beta_0(\psi,\phi)
\end{equation*}
\begin{equation}\label{3.6}
    \leq c{\mathscr N}(\psi^\prime)+b\alpha(\psi^\prime,\phi^\prime)-\dot\beta_0(\psi^\prime,\phi^\prime)
    \leq c{\mathscr N}+b\alpha -\dot\beta_0(\psi^\prime,\phi^\prime)
\end{equation}

where to get the last inequality we used (\ref{3.4}).

It follows from (\ref{3.6}) that
$$\dot\beta_0(\psi,\phi)\geq \dot\beta_0(\psi^\prime,\phi^\prime).
$$

To get the last statement of the theorem we note that if
$\dot\beta_0(\psi,\phi)=\dot\beta_0(\psi^\prime,\phi^\prime)$
then there are equalities in
(\ref{3.6}) instead of the inequalities which is only
possible if
${\mathscr N}(\psi^\prime)={\mathscr N}$ and $\alpha(\psi^\prime,\phi^\prime)=\alpha$.
\rule{2mm}{2mm}

\begin{remark}\label{r3.1}
There is a more restricted definition of locally most powerful tests in \cite{berk}, \cite{roters}, \cite{schmitz} (among others),  where the derivative of the power function is maximized within the class of all tests with a given value of the type I error probability ($\alpha(\psi^\prime,\phi^\prime)=\alpha$ instead of $\alpha(\psi^\prime,\phi^\prime)\leq\alpha$ in (\ref{3.4})). It is obvious that, with this modification, the assertion of Theorem  \ref{t3.1}  is also valid if the conditions of Theorem \ref{t3.1} hold with $b<0$.

If the conditions of Theorem \ref{t3.1} hold with $b=0$, then  for any test $(\psi^\prime,\phi^\prime)\in\Delta$ satisfying
\begin{equation}\label{3.7}
{\mathscr N}(\psi^\prime)\leq{\mathscr N}
 \end{equation}
 it holds
\begin{equation}\label{3.8}
\dot\beta_0(\psi,\phi)\geq  \dot\beta_0(\psi^\prime,\phi^\prime). \end{equation}

 The inequality
in (\ref{3.8}) is strict if the inequality in
(\ref{3.7}) is strict.
\end{remark}

For any stopping rule $\psi=(\psi_1,\psi_2,\dots)$ let us denote
\begin{equation*}
s_n^\psi=(1-\psi_1)\dots (1-\psi_{n-1})\psi_n\quad\mbox{and}\quad t_n^\psi=(1-\psi_1)\dots (1-\psi_{n-1}),
\end{equation*}
for any $n=1,2,\dots$ ($t_1^\psi\equiv 1$ by definition).

Let $I_A$ be the indicator function of the event $A$.

The following theorem, in a rather standard way (see, for example, \cite{berk} or \cite{schmitz}), lets us find optimal decision rules for any given stopping rule $\psi$.

\begin{theorem} \label{t3.2} Let Assumption 1 be fulfilled.
For any $b\in\mathbb R$  and for any sequential test $(\psi,\phi)$
such that $E_0\tau_\psi<\infty$
\begin{equation}\label{3.9}
b\alpha(\phi,\psi)- \dot\beta_0(\phi, \psi)\geq
\sum_{n=1}^\infty\int s_n^\psi \min \{0,b f_0^n-\dot f_0^n\} d\mu^n
\end{equation}
with an equality if and only if
\begin{equation}\label{3.10}
I_{\{bf_0^n<\dot f_0^n\}}\leq\phi_n\leq I_{\{bf_0^n\leq\dot f_0^n\}}
\hspace{12 pt}
\end{equation}
$\mu^n$-almost everywhere on $S_n^\psi=\left\{(x_1,\dots,x_n): s_n^\psi(x_1,\dots,x_n)>0\right\}$ for any $n=1,2,\dots$.
\end{theorem}

The proof of Theorem \ref{t3.2} can be found in Appendix.

Let us denote
$$L(\psi)=L(\psi;b,c)=\inf_\phi L(\psi,\phi;b,c).$$

\begin{corollary}\label{c3.1} Under conditions of Theorem \ref{t3.2}, if $E_0\tau_\psi<\infty$, then
\begin{equation}\label{3.11}
L(\psi)=\sum_{n=1}^\infty \int s_n^\psi(cnf^n+l_n) d\mu^n,
\end{equation}
where, by definition,
\begin{equation*}\label{21}l_n=\min\{0,bf_0^n-\dot f_0^n\}.
\end{equation*}
\end{corollary}
{\em Proof.}
This follows from Theorem \ref{t3.2} by (\ref{reduction5}), in view of
(\ref{1.2}).
\rule{2mm}{2mm}\\

By Theorem \ref{t3.2}, the problem of minimization of $L(\psi,\phi;b,c)$
is reduced now to the problem of minimization of $L(\psi;b,c)$, that is,
to an optimal stopping problem. Indeed, if there is a $\psi$ such that $E_0\tau_\psi<\infty$ and such that
$$L(\psi;b,c)=\inf_{\psi^\prime}L(\psi^\prime;b,c),$$
then, adding to $\psi$ any decision rule $\phi$ satisfying
(\ref{3.10}), by Theorem \ref{t3.2} we have that for any sequential
test $(\psi^\prime,\phi^\prime)$:
\begin{equation*}
L(\psi,\phi;b,c)=L(\psi;b,c)\leq L(\psi^\prime;b,c)\leq
L(\psi^\prime,\phi^\prime;b,c).
\end{equation*}
In particular, in this way we obtain tests $(\psi,\phi)$ satisfying
(\ref{3.2}), which is crucial for solving the original conditional
problem (see Theorem \ref{t3.1}).

\section{Optimal Stopping Rules}
In this section, we characterize the structure of  stopping rules minimizing $L(\psi)$, first in the class of truncated stopping rules, then in some natural classes of non-truncated stopping rules.

We suppose, throughout this Section, that Assumption 1 is fulfilled.
\subsection{\normalsize Optimal Truncated Stopping Rules}

Here we solve the problem of
minimization of $L(\psi)$ in the class of truncated stopping rules,
that is, in the class ${\mathscr F}^N$, $N\geq 1$, of stopping rules $\psi$
such that
\begin{equation}\label{optimal5}(1-\psi_1)(1-\psi_2)\dots(1-\psi_N)\equiv 0.\end{equation}

For any stopping rule $\psi\in{\mathscr F}^N$ let
\begin{equation}\label{4.2}
L_N(\psi)=L_N(\psi;b,c)=\sum_{n=1}^{N-1}\int
s_n^\psi(cnf^n+l_n)d\mu^n
+\int t_N^\psi\left(cNf^N+l_N\right)d\mu^N
\end{equation}
(it is easy to see that, by virtue of (\ref{optimal5}), $L_N(\psi)$ coincides with the right-hand side of (\ref{3.11})).

Let us define
$V_N^N\equiv l_N$, and
recursively for any $n=N-1,N-2,\dots 1$
\begin{equation}\label{4.3}
V_n^N=\min\{l_n,cf^n+R_n^N\},
\end{equation}
where
\begin{equation*}
R_n^N=R_n^N(x_1,\dots,x_n)=\int V_{n+1}^N(x_1,\dots,x_{n+1})d\mu(x_{n+1}).
\end{equation*}
Let, finally,  for any $k=1,\dots, N$
\begin{equation*}
Q_k^N(\psi)=\sum_{n=1}^{k-1}\int
s_n^\psi(cnf^n+l_n)d\mu^n
+\int t_k^\psi\left(ckf^k+V_k^N\right)d\mu^k
\end{equation*}

\begin{theorem}\label{t4.1} Let $\psi\in{\mathscr F}^N$ be any (truncated) stopping rule. Then for any
$1\leq k\leq N$
\begin{equation}\label{4.4}
L_N(\psi)\geq Q_k^N(\psi)
\end{equation}

The lower bound in (\ref{4.4}) is attained if and only if
\begin{equation}\label{4.5}
I_{\{l_{n}< cf^{n}+R_n^N\}}\leq\psi_{n}\leq I_{\{l_{n}\leq cf^{n}+R_n^N
\}}
\end{equation}
$\mu^n$-almost everywhere on $$T_n^\psi=\{(x_1,\dots, x_n):t_n^\psi(x_1,\dots, x_n)>0\}$$ for any $n=k,k+1,\dots, N-1$.
\end{theorem}

The proof of Theorem \ref{t4.1} is laid down in Appendix.

From Theorem \ref{t4.1} we easily have  the following optimality result.
\begin{corollary}\label{c4.1}
 For any $\psi\in {\mathscr F}^N$
\begin{equation}\label{4.6}
L_N(\psi)\geq c+R_0^N,
\end{equation}
where
$$
R_0^N=\int V_1^N(x_1)d\mu(x_1).
$$
There is an equality in (\ref{4.6})
if and only if $\psi_n$ satisfy (\ref{4.5}) $\mu^n$-almost everywhere on $T_n^\psi$, for any $n=1,2,\dots, N-1$.
\end{corollary}

{\em Proof.}
It is straightforward because $$Q_1^N(\psi)=\int t_1^\psi\left(cf^1+V_1^N\right)d\mu=c+\int V_1^N d\mu=c+R_0^N.$$
\rule{2mm}{2mm}

\subsection{Optimal Non-Truncated Stopping Rules}

In this section we characterize the structure of general
sequential tests minimizing $L(\psi)=L(\psi;b,c)$ (see (\ref{reduction5})).

 Let us define for any stopping rule $\psi$, and for any natural $N\geq 1$,
\begin{equation*}L_N(\psi)=L_N(\psi;b,c)=L(\psi^N;b,c),
\end{equation*}
where $\psi^N=(\psi_1,\psi_2,\dots,\psi_{N-1},1,\dots)$ is the rule $\psi$ truncated at $N$.

By (\ref{4.2}),
\begin{equation*}
L_N(\psi)=\sum_{n=1}^{N-1}\int
s_n^\psi(cnf^n+l_n)d\mu^n
+\int t_N^\psi\left(cNf^N+l_N\right)d\mu^{N}.
\end{equation*}

Because $\psi^N$ is truncated, the results of the preceding section apply, in particular, Theorem \ref{t4.1}.
The idea of the following construction is to pass to the limit, as $N\to\infty$, in (\ref{4.4}), in order to get some lower bound for the "risk" $L(\psi)$, and corresponding conditions under which the lower bound is attained.

First of all, let us show that the right-hand side of (\ref{4.4}) has a limit, as $N\to\infty$, for any $k=1,2,\dots$.
This is basically due to the following

\begin{lemma}\label{l4.1} For any $n\geq 1$ and for any  $N\geq n$
\begin{equation}\label{4.7}
l_n\geq V_n^N\geq V_n^{N+1}.
\end{equation}
\end{lemma}

The first inequality in (\ref{4.7}) is due to (\ref{4.3}). The proof of the second is identical to the proof of Lemma 3.3 \cite{novikov09b}.

It follows from Lemma \ref{l4.1} that for any fixed $n\geq 1$ the
sequence $V_n^N$, $N=1,2,\dots$, is non-increasing. So, there exists
\begin{equation}\label{4.8}V_n= \lim_{N\to\infty}V_n^N
\end{equation}
such that $V_n\leq l_n$ for all $n=1,2,\dots$.

Because of this,  the right-hand side of (\ref{4.4}), by the Lebesgue's monotone convergence theorem converges to
\begin{equation}\label{4.9}
Q_k(\psi)=\sum_{n=1}^{k-1}\int
s_n^\psi(cnf^n+l_n)d\mu^n
+\int t_k^\psi\left(ckf^k+V_k\right)d\mu^k
\end{equation}
for any $k=1,2,\dots$.
By the same reason, passing to the limit on both sides of (\ref{4.3}) is possible, which gives us
\begin{equation*}V_n=\min\{l_n,cf^n+R_n\},\end{equation*}
where
\begin{equation*}R_n=R_n(x_1,\dots,x_n)=\int V_{n+1}(x_1,\dots,x_{n+1})d\mu(x_{n+1}),\end{equation*}
for any $n=1,2,\dots$.

At last, to be able to pass to the limit on the left-hand side of (\ref{4.4}), we need that
$L_N(\psi)\to L(\psi)$, as $N\to\infty$, at least for some class of stopping rules $\psi$.
Let $\mathscr F$ be a class of stopping rules  such that for every $\psi\in\mathscr F$
it holds
\begin{equation}\label{4.12}
E_0 \tau_\psi<\infty,\;E\tau_\psi<\infty,\quad\mbox{and}\;\lim_{N\to\infty}L_N(\psi;b,c)=L(\psi;b,c)
\end{equation}
for all $b\in\mathbb R$ and $c>0$ (the first condition in (\ref{4.12}) is needed in order that (\ref{3.11}) be valid, the second one guarantees that $L(\psi;b,c)<\infty$).

Now passing to the limit on both sides of  (\ref{4.4}), as $N\to\infty$,  is possible for all $\psi\in\mathscr F$, so we get

\begin{lemma}\label{l4.2} For any stopping rule $\psi\in\mathscr F$ and for
any $k\geq 1$
\begin{equation*}
L(\psi)\geq Q_k(\psi),
\end{equation*}
where $Q_k(\psi)$ is defined by (\ref{4.9}), being $V_n$ defined, for any  $n=1,2,\dots$, by
(\ref{4.8}).

In particular,  for any stopping rule $\psi\in\mathscr F$
\begin{equation}\label{4.13}
L(\psi)\geq c+R_0.
\end{equation}
 \end{lemma}

The following lemma shows that the lower bound in  (\ref{4.13}) is, in fact, the infimum value of the left-hand side of (\ref{4.13}).

\begin{lemma}\label{l4.3} Let $\mathscr G\subset\mathscr F$ be any subclass of stopping rules, such that
\begin{equation}\label{4.14}
    \bigcup_{N=1}^\infty\mathscr F^N\subset \mathscr G.
\end{equation}
Then
\begin{equation*}
\inf_{\psi\in \mathscr G}L(\psi)=c+R_0.
\end{equation*}
\end{lemma}
{\em Proof.}
If $R_0>-\infty$, then the proof is conducted in the same way as the proof of Lemma 3.5 in \cite{novikov09b}.

If $R_0=-\infty$, it follows from $\lim_{N\to\infty}R_0^N=R_0$ that for any $k>-\infty$ there is $N$ such that $R_0^N\leq k$. Thus, because of (\ref{4.14}), $\inf_{\psi\in\mathscr G}L(\psi)\leq k$.
Because $k>-\infty$ is arbitrary,  $\inf_{\psi\in\mathscr G}L(\psi)=-\infty$ follows.
\rule{2mm}{2mm}

\begin{theorem}\label{t4.2} Let Assumption 2 be fulfilled and let $\mathscr G$ be any class of stopping rules satisfying
the conditions of Lemma \ref{l4.3} and such that
\begin{equation}\label{4.15}
\inf_{\psi^\prime\in \mathscr G}L(\psi^\prime)>-\infty.
\end{equation}

If there exists $\psi$ in ${\mathscr G}$
such that
\begin{equation}\label{4.16}
L(\psi)=\inf_{\psi^\prime\in \mathscr G}L(\psi^\prime),
\end{equation}
then
\begin{eqnarray}\label{4.17}
I_{\{l_{n}< cf^{n}+R_n \}}\leq \psi_{n}\leq I_{\{l_{n}\leq cf^{n}+R_n
\}}
\end{eqnarray}
$\mu^n$-almost everywhere on $T_n^\psi$, for any $n=1,2,\dots$, and
\begin{equation}\label{4.18}
\int t_n^\psi (V_n-l_n)d\mu^n\to 0, \mbox{ as } n\to\infty.
\end{equation}

On the other hand, if  $\psi$ satisfies (\ref{4.17}) $\mu^n$-almost everywhere on $T_n^\psi$, for any $n=1,2,\dots$, and  satisfies (\ref{4.18}), and if $\psi\in\mathscr G$, then it satisfies (\ref{4.16}) as well.
\end{theorem}

The proof of Theorem \ref{t4.2} can be found in Appendix.

\begin{remark} Generally speaking, (\ref{4.15}) can be violated.   Let us see the following example, in which
$$
\inf_{\psi\in \mathscr G}L(\psi;b,c)=-\infty
$$
for all $b\in\mathbb R$ and $c>0$.

Suppose that  $X_1, X_2,\dots$ are independent and that $X_n$ is normally distributed with mean $n\theta$ and unit variance ($X_n\sim\mathscr N(n\theta, 1)$), $n=1,2\dots$.  Suppose also that $H_0:\theta=0$ and $H_1:\theta>0$.

Let $\psi^N$ be a fixed sample size-stopping rule taking $N$ observations ($\psi_1^N=\hdots=\psi_{N-1}^N=0,\psi_N^N=1$), $N=1,2,\dots$. Then it is easy to see that
$$L(\psi^N;b,c)=cN+b\,\Phi(-b/\sigma_N)-\frac{\sigma_N}{\sqrt{2\pi}}\exp\{-b^2/2\sigma_N^2\},$$
where $\sigma_N^2=\sum_{n=1}^N n^2\sim N^3/3$, thus,
$L(\psi^N;b,c)\to-\infty$ as $N\to\infty$, for any $b\in\mathbb R$ and $c>0$.
\end{remark}

With respect to the  property of (\ref{4.15}), any hypothesis testing problem exhibits the following remarkable behavior.

\begin{theorem}\label{t4.3} For any  family $\{f_\theta^n,\,\theta\in\Theta,\,n=1,2,\dots\}$ of the finite-dimensional distributions of the process $X_1,X_2,\dots$, satisfying Assumption 2, either
\begin{equation}\label{4.19}
\inf_{\psi\in \mathscr G}L(\psi;b,c)>-\infty
\end{equation}
for all $b\in\mathbb R$ and $c>0$, or
\begin{equation*}
\inf_{\psi\in \mathscr G}L(\psi;b,c)=-\infty
\end{equation*}
for all $b\in\mathbb R$ and $c>0$
\end{theorem}

The proof of Theorem \ref{t4.3} is laid down in Appendix.

In view of Theorem \ref{t4.3} the following definition is justified. We call a hypothesis testing problem {\em finite} if (\ref{4.19}) is fulfilled for all $b\in\mathbb R$ and $c>0$. For hypothesis testing problems which are not finite, we do not have any other recommendation than minimization of the Lagrange multiplier function $L_N(\psi;b,c)$, for some $b\in \mathbb R$ and $c>0$, in the class $\mathscr F^N$ of truncated stopping rules using Corollary \ref{c4.1}. For finite problems, we may hope to find optimal non-truncated stopping rules using Theorem \ref{t4.2} (see Section \ref{s5} below).

There is a way to make the sufficient condition of optimality in Theorem \ref{t4.2}  more practical, supposing that, additionally to Assumption 2,  Assumption 3 holds. Namely, it can be shown that in this case $\mathscr G_1=\{\psi:E\tau_\psi<\infty, E_0\tau_\psi<\infty\}\subset \mathscr F$ satisfies the conditions of Lemma \ref{l4.3}. Therefore, by Theorem \ref{t4.2}, any  $\psi$ satisfying (\ref{4.17}) and (\ref{4.18}) will be optimal if $\psi\in \mathscr G_1$. We formalize this in the following
\begin{lemma}\label{l4.4}
Let us suppose that Assumptions 2 and 3 are satisfied and that the hypothesis testing problem is finite. Then \begin{equation}\label{4.21}\bigcup_{N=1}^\infty \mathscr F^N\subset\mathscr G_1\subset \mathscr F.\end{equation}
\end{lemma}

The proof of Lemma \ref{l4.4} is laid down in Appendix.

\section{Applications to conditional problems}\label{s5}


For any $c>0$ and $b\in\mathbb R$ let us call a sequential test $(\psi,\phi)$ $(b,c)$-{\em generated} if the following conditions are fulfilled:

\begin{equation}\label{5.1}
 I_{\{l_n<cf^n+\int V_{n+1}d\mu(x_{n+1})\}}\leq \psi_n\leq I_{\{l_n\leq cf^n+\int V_{n+1}d\mu(x_{n+1})\}}
 \end{equation}
 $\mu^n$-almost everywhere on $T_n^\psi$, $n=1,2,\dots$,
 and
  \begin{equation}\label{5.2}
     I_{\{bf_0^n<\dot f_0^n\}}\leq \phi_n\leq I_{\{bf_0^n\leq\dot f_0^n\}}
  \end{equation}
  $\mu^n$-almost everywhere on $S_n^\psi$, $n=1,2,\dots$, where

  \begin{equation}\label{5.3}
    l_n=\min\{0,bf_0^n-\dot f_0^n\}
 \end{equation}
  \begin{equation*}
    V_n=\lim_{N\to\infty}V_n^N,
 \end{equation*}for $n=1,2,\dots$,
 being $V_n^N$ defined recursively, for $n=N-1,N-2,\dots, 1$, by
  \begin{equation}\label{5.5}
    V_n^N=\min\{l_n,cf^n+\int V_{n+1}^Nd\mu(x_{n+1})\},
 \end{equation}
 starting from $V_N^N\equiv l_N$, $N=1,2,\dots$.

Let us call a $(b,c)$-generated test $(\psi,\phi)$ {\em regular} if
 \begin{equation*}
    \int t_n^\psi (V_n-l_n)d\mu^n\to 0,\quad n\to\infty,
 \end{equation*}
  holds.

The following Theorem is an immediate consequence of Theorems \ref{t3.1}, \ref{t3.2} and \ref{t4.2}.
\begin{theorem}\label{t5.1}
 Let the conditions of Theorem \ref{t4.2} be satisfied. Let $c>0$ and $b>0$ be arbitrary constants, and
 let  $(\psi,\phi)$ be any $(b,c)$-generated regular test, such that
 $\psi\in \mathscr G$.

  Then the sequential test $(\psi,\phi)$ is locally most powerful for testing $H_0:\theta=\theta_0$ vs $H_1:\theta>\theta_0$ in the following sense.

  For any $(\psi^\prime,\phi^\prime)$, with $\psi^\prime\in \mathscr G$, such that
  \begin{equation}\label{5.7}
   E\tau_{\psi^\prime}\leq E\tau_\psi,\quad\mbox{and}\quad\alpha(\psi^\prime,\phi^\prime)\leq\alpha(\psi,\phi),
  \end{equation}
  it holds
  \begin{equation}\label{5.8}
    \dot \beta_0(\psi,\phi)\geq \dot \beta_0(\psi^\prime,\phi^\prime).
  \end{equation}
  The inequality in (\ref{5.8}) is strict if at least one of the inequalities in (\ref{5.7}) is strict.

  If there are equalities in all of the inequalities in (\ref{5.7}) and (\ref{5.8}), then $(\psi^\prime,\phi^\prime)$ is a $(b,c)$-generated regular test as well.
 \end{theorem}

 \begin{remark}\label{r5.1}
 Under the conditions of  Theorem \ref{t5.1}, if  $(\psi,\phi)$ is a $(b,c)$-generated regular test with $b<0$, and $\psi\in\mathscr G$, then it is locally most powerful in the class of all sequential tests of the (exact) size $\alpha(\psi,\phi)$ (with an equality instead of the second inequality in (\ref{5.7}), see Remark \ref{r3.1}).

 Similarly, if, under the same conditions,  if $(\psi,\phi)$ is $(b,c)$-generated regular test with  $b=0$, and $\psi\in\mathscr G$, then it is is locally most powerful in the sense that for any $(\psi^\prime,\phi^\prime)$, with $\psi^\prime\in\mathscr G$, such that
 $$
 E \tau_{\psi^\prime}\leq E \tau_{\psi}
 $$
 it holds
 $$
 \dot \beta_0(\psi,\phi)\geq\dot \beta_0(\psi^\prime,\phi^\prime),
 $$
 irrespective of the corresponding type I error probabilities, with the respective modification of Theorem \ref{t5.1}.
 \end{remark}
 It is interesting to note that if $(\psi,\phi)$ is $(b,c)$-generated regular test with $b<0$, then the test $(\psi,\bar\phi)$, where, by definition, $\bar\phi=(1-\phi_1,1-\phi_2,\dots)$, is locally most powerful, in the sense of Theorem \ref{t5.1}, for testing $H_0:\theta=\theta_0$ vs $H_1:\theta<\theta_0$.
 To make this formal, we need some additional results.

 Let for any $b\in\mathbb R$ and $c>0$
 $$
 \bar L(\psi)= \bar L(\psi;b,c)=\inf_{\phi} (cE\tau_\psi+b\alpha(\psi,\phi)+\dot\beta_0(\psi,\phi)),
 $$
 and $\bar L_N(\psi;b,c)=\bar L(\psi^N;b,c)$, where $\psi^N=(\psi_1,\dots,\psi_{N-1},1,\dots)$.
 \begin{theorem}\label{t5.2} Suppose that Assumption 2 is fulfilled.

  Let $c>0$ and $b\in\mathbb R$ be arbitrary constants.

 Let $\mathscr G$ be some class of tests  such that
 $$
 \bigcup_{N=1}^\infty \mathscr F^N\subset \mathscr G
 $$
 and such that $L_N(\psi;b,c)\to L(\psi;b,c)$ and $\bar L_N(\psi;-b,c)\to \bar L(\psi;-b,c)$, as $N\to\infty$, for all $\psi\in \mathscr G$.

 Then for a stopping rule $\psi\in\mathscr G$
 \begin{equation*} L(\psi;b,c)=\inf_{\psi^\prime\in\mathscr G}  L(\psi^\prime;b,c)
 \end{equation*}
 if and only if
  \begin{equation*} \bar L(\psi;-b,c)=\inf_{\psi^\prime\in\mathscr G}  \bar L(\psi^\prime;-b,c).
 \end{equation*}
 \end{theorem}
The proof of Theorem \ref{t5.2} can be found in Appendix.

Using Theorem \ref{t3.1} and Theorem \ref{t3.2}, we get from Theorem \ref{t5.2} the following
 \begin{theorem}\label{t5.3}
 Let the conditions of Theorem \ref{t5.2} be satisfied with some $c>0$ and $b<0$.
 Let $(\psi,\phi)$ be any $(b,c)$-generated regular test, such that
 $\psi\in \mathscr G$, and let the problem of testing $H_0:\theta=\theta_0$ vs $H_1:\theta>\theta_0$ be finite.

 Then the sequential test $(\psi,\bar\phi)$ is locally most powerful for testing $H_0:\theta=\theta_0$ vs $H_1:\theta<\theta_0$ in the following sense.

  For any $(\psi^\prime,\phi^\prime)$, with $\psi^\prime\in \mathscr G$, such that
  \begin{equation}\label{5.11}
   E\tau_{\psi^\prime}\leq E\tau_\psi\quad\mbox{and}\quad\alpha(\psi^\prime,\phi^\prime)\leq\alpha(\psi,\bar\phi),
  \end{equation}
  it holds
  \begin{equation}\label{5.12}
    \dot \beta_0(\psi,\bar\phi)\geq \dot \beta_0(\psi^\prime,\phi^\prime).
  \end{equation}
  The inequality in (\ref{5.12}) is strict if at least one of the inequalities in (\ref{5.11}) is strict.

  If there are equalities in all of the inequalities in (\ref{5.11}) and (\ref{5.12}), then $(\psi^\prime,\bar\phi^\prime)$ is a regular $(b,c)$-generated test as well.\end{theorem}

In the rest of this section, we will apply the results of Theorems \ref{t5.1} and \ref{t5.3} to the case of i.i.d. observations considered in  \cite{berk}. Obviously, the conditions of \cite{berk} imply that our Assumptions 1 to 3 are fulfilled, thus, we can make use of all our results above.
In this case $f_\theta^n=f_\theta^n(x_1,\dots,x_n)=\prod_{i=1}^n f_\theta(x_i)$, where $f_\theta$ is the marginal density with respect to $\mu$, and $$\left.\dot f_0^n(x_1,\dots,x_n)=(f_{\theta}^n(x_1,\dots,x_n))_\theta^\prime\right|_{\theta=\theta_0}.$$

As in \cite{berk}, we are using $E\tau_\psi=E_{\theta_0}\tau_\psi$ in the conditional minimization problems.

Let us see first, how the structure of $(b,c)$-generated tests transforms in this case (see (\ref{5.3}) - (\ref{5.5})).

It is immediate that $l_n=\min\{0,b-z_n\}f_0^n$, where \begin{equation}\label{5.13}z_n=\sum_{i=1}^n \frac{( f_{\theta}(x_i))_\theta^\prime|_{\theta=\theta_0}}{f_{\theta_0}(x_i)}\end{equation}
(we use here the fact that $\dot f_0^n=0$ $\mu^n$-almost everywhere on $\{f_{\theta_0}^n=0\}$, which easily follows from (\ref{2.1})). The definition of $z_n$ in case $f_{\theta_0}(x_i)=0$ does not matter, because in this case $f_0^n(x_1,\dots,x_n)=\prod_{i=1}^n f_{\theta_0}(x_i)=0$.

Let the $i$-th summand on the right-hand side of (\ref{5.13}) be denoted as $r_i=r(x_i)$.

Let us define $g(z)\equiv \min\{0,-z\}$, $z\in\mathbb R$. Let further $\rho_c^0(z)=g(z)$, $z\in\mathbb R$, and for any $n=1,2,\dots$, recursively,
\begin{equation}\label{5.14}
\rho_c^n(z)=\min\{g(z),c+\int\rho_c^{n-1}(z+r(x))f_{\theta_0}(x)d\mu(x)\},
\end{equation}
$z\in\mathbb R$.

It is easy to see, by induction, that $V_n^N=\rho_c^{N-n}(z_n-b)f_0^n$, for any $n=N,N-1,\dots 1$, and for any $N=1,2,\dots$.
It is also easy to see, by induction, that $\rho_c^n(z)\geq \rho_c^{n+1}(z)$ for any $z\in \mathbb R$ and for any $n=0,1,\dots$.
Thus, there exists $\rho_c(z)=\lim_{n\to\infty}\rho_c^n(z)$, $z\in \mathbb R$. Below, we will prove that $\rho_c(z)$ is finite for any $z\in\mathbb R$.

Therefore, $V_n=\lim_{N\to\infty}V_n^N=\lim_{N\to\infty}\rho_c^{N-n}(z_n-b)f_0^n=\rho_c(z_n-b)f_0^n$. In particular,
$V_1=\rho_c(z_1-b)f_0^1=\rho_c(r_1-b)f_{\theta_0}(x_1)$, thus, \begin{equation*} R_0=\int V_1(x_1)d\mu(x_1)=\int \rho_c(r(x)-b)f_{\theta_0}(x)d\mu(x).\end{equation*}

Further, passing to the limit in (\ref{5.14}), as $n\to\infty$, we have
\begin{equation}\label{5.16}
\rho_c(z)=\min\{g(z),c+\int\rho_c(z+r(x))f_{\theta_0}(x)d\mu(x)\},
\end{equation}

The inequality $l_n\leq cf_0^n+\int V_{n+1}d\mu(x_{n+1})$ in (\ref{5.1}) is equivalent now to
$$\rho_c(z_n-b)\leq c+\int \rho_c(z_n-b+r(x))f_{\theta_0}(x)d\mu(x) $$
on $\{f_0^n>0\}$. Respectively, the inequality $bf_0^n\leq \dot f_0^n$ in (\ref{5.2}) is equivalent to $b\leq z_n$ on $\{f_0^n>0\}$.

It follows that a sequential test $(\psi,\phi)$ is $(b,c)$-generated if and only if

\begin{equation}\label{5.17}
 I_{\{g(z_n-b)<c+\int \rho_c(z_n-b+r(x))f_{\theta_0}(x)d\mu(x)\}}\leq \psi_n\leq I_{\{g(z_n-b)\leq c+\int \rho_c(z_n-b+r(x))f_{\theta_0}(x)d\mu(x)\}}
 \end{equation}
 $\mu^n$-almost everywhere on $T_n^\psi\cap\{f_0^n>0\}$, $n=1,2,\dots$,
 and
  \begin{equation*}
     I_{\{b<z_n\}}\leq \phi_n\leq I_{\{b\leq z_n\}}
  \end{equation*}
  $\mu^n$-almost everywhere on $S_n^\psi\cap\{f_0^n>0\}$, $n=1,2,\dots$.

  Respectively, a $(b,c)$-generated test $(\psi,\phi)$ is regular if
  \begin{equation}\label{5.19}
    \int t_n^\psi (\rho_c(z_n-b)-g(z_n-b))f_0^nd\mu^n\to 0,\quad n\to\infty.
 \end{equation}

 The plan of the rest of this section is as follows.
 Let $$h_c(z)=\int \rho_c(z+r(x))f_{\theta_0}(x)d\mu(x)$$
 (see \ref{5.16}). Then, if $c+h_c(0)\leq 0$, it can be shown that
  there exist $A_c\leq 0$ and $B_c\geq 0$ such that \begin{equation*} g(z)=c+h_c(z)\end{equation*} for $z=A_c$ and $z=B_c$, and that
 the inequality $g(z)>c+h_c(z)$ is equivalent to $z\in(A_c,B_c)$.

 Thus, it will follow that a test $(\psi,\phi)$ is $(b,c)$-generated (supposing that $c+h_c(0)\leq 0$) if and only if
 \begin{equation}\label{5.21}
 I_{\{z_n\not\in [b+A_c,\,b+B_c]\}}\leq \psi_n\leq I_{\{z_n\not\in (b+A_c,\,b+B_c)\}}
 \end{equation}
 $\mu^n$-almost everywhere on $T_n^\psi\cap\{f_0^n>0\}$, $n=1,2,\dots$,
 and
  \begin{equation}\label{5.22}
   I_{\{z_n> b\}}  \leq \phi_n\leq I_{\{z_n\geq b\}}
  \end{equation}
  $\mu^n$-almost everywhere on $S_n^\psi\cap\{f_0^n>0\}$, $n=1,2,\dots$.

  If $c+h_c(0)>0$, it will follow that a test $(\psi,\phi)$ is $(b,c)$-generated if and only if
  $\psi_n=1$
 $\mu^n$-almost everywhere on $T_n^\psi\cap\{f_0^n>0\}$, $n=1,2,\dots$,
 and (\ref{5.22}) is satisfied
  $\mu^n$-almost everywhere on $S_n^\psi\cap\{f_0^n>0\}$, $n=1,2,\dots$.

   Therefore,  in this particular case from Theorem \ref{t5.1} we will have
 \begin{theorem}\label{t5.4}
 Let Assumptions 1 -- 4 of \cite{berk} be fulfilled. Let $c>0$  and  $b>0$ be any constants, and
 let  $(\psi,\phi)$ be any $(b,c)$-generated test.

  Then $E_0\tau_\psi<\infty$, and the sequential test $(\psi,\phi)$ is locally most powerful for testing $H_0:\theta=\theta_0$ vs $H_1:\theta>\theta_0$ in the following  sense.
   For any $(\psi^\prime,\phi^\prime)$ such that
  \begin{equation}\label{5.23}
   E_0 \tau_{\psi^\prime}\leq E_0 \tau_\psi,\quad\mbox{and}\quad\alpha(\psi^\prime,\phi^\prime)\leq\alpha(\psi,\phi),
  \end{equation}
  it holds
  \begin{equation}\label{5.24}
    \dot \beta_0(\psi,\phi)\geq \dot \beta_0(\psi^\prime,\phi^\prime).
  \end{equation}

  The inequality in (\ref{5.24}) is strict if at least one of the inequalities in (\ref{5.23}) is strict.

  If there are equalities in all of the inequalities in (\ref{5.23}) and (\ref{5.24}), then $(\psi^\prime,\phi^\prime)$ is a $(b,c)$-generated test as well.
 \end{theorem}

  The proof of Theorem \ref{t5.4} can be found in Appendix.

 Analogously, from Theorem \ref{t5.3} we obtain
 \begin{theorem}\label{t5.5}
 Let Assumptions 1 -- 4 of \cite{berk} be fulfilled. Let $c>0$ and $b<0$ be arbitrary constants, and
 let  $(\psi,\phi)$ be any $(b,c)$-generated test.

  Then $E_0\tau_\psi<\infty$, and the sequential test $(\psi,\bar\phi)$ is locally most powerful for testing $H_0:\theta=\theta_0$ vs $H_1:\theta<\theta_0$ in the following  sense.
   For any $(\psi^\prime,\phi^\prime)$ such that
  \begin{equation}\label{5.25}
   E_0\tau_{\psi^\prime}\leq E_0\tau_\psi,\quad\mbox{and}\quad\alpha(\psi^\prime,\phi^\prime)\leq\alpha(\psi,\bar\phi),
  \end{equation}
  it holds
  \begin{equation}\label{5.26}
    \dot \beta_0(\psi,\bar\phi)\geq \dot \beta_0(\psi^\prime,\phi^\prime).
  \end{equation}

  The inequality in (\ref{5.26}) is strict if at least one of the inequalities in (\ref{5.25}) is strict.

  If there are equalities in all of the inequalities in (\ref{5.25}) and (\ref{5.26}), then $(\psi^\prime,\bar\phi^\prime)$ is a $(b,c)$-generated test as well.
 \end{theorem}

 \begin{remark}\label{r5.2}
 It can be shown (very much like in the proof of Theorem 6 in \cite{novikov08}, see also \cite{berk} for the non-randomized case) that for any $-\infty<A<B<\infty$ any sequential test $(\psi,\phi)$ with
 \begin{equation}\label{5.27}
 I_{\{z_n\not\in[A,B]\}}\leq\psi_n\leq I_{\{z_n\not\in(A,B)\}}\end{equation}
 and
 \begin{equation}\label{5.28}
 \phi_n=I_{\{z_n\geq B\}},
 \end{equation}
 $n=1,2,\dots$, is $(b,c)$-generated for some $c>0$ and $b\in\mathbb R$.

 \cite{roters} (see Remark i) on page 182) notes that, generally speaking, a test of type (\ref{5.27})-(\ref{5.28}) is not locally most powerful (in the sense of our Theorem \ref{t5.4}) and gives an example of a test $(\psi,\phi)$ of type (\ref{5.27})-(\ref{5.28}), for which there exists another test $(\psi^\prime,\phi^\prime)$ such that $E_0 \tau_\psi=E_0 \tau_{\psi^\prime}=k^2$, and $\alpha(\psi^\prime,\phi^\prime)=0.5<\alpha(\phi,\psi)=0.8$ and  $\dot\beta(\psi^\prime,\phi^\prime)=k>\dot\beta(\phi,\psi)=0.8k$, where $k$ is any natural number.  It follows from Theorem \ref{t5.4} that the only way this can happen is that $(\psi,\phi)$ is $(b,c)$-generated with $b<0$.
 Thus, it follows from Theorem \ref{t5.5} that $(\psi,\bar\phi)$ is locally most powerful for testing $\theta=\theta_0$ vs $\theta<\theta_0$ at level $\alpha(\psi,\bar \phi)=1-\alpha(\psi,\phi)=0.2$ of significance.

 It is interesting to note that $(\psi^\prime,\phi^\prime)$ in his example is $(b,c)$-generated with $b=0$, so it is locally most powerful among all sequential tests with the same, or lesser, average sample number, irrespective of their $\alpha$-level (see Remark \ref{r5.1}).
 \end{remark}
 \begin{remark}
 It is easy to see that if the distribution of $r(X_1)$ is symmetric under $H_0$ (as, for example, in the case of normal family $\mathscr N(\theta,\sigma^2),\,\theta\in\mathbb R$, $\sigma^2>0$), then $A_c=-B_c$, i.e. the continuation region of any $(b,c)$-generated test is symmetric with respect to $b$. In such a case, it follows from Theorems {\ref{t5.4}} and \ref{t5.5} that any test $(\psi,\phi)$ with $$I_{\{z_n\not\in[b-B_c,b+B_c]\}}\leq\psi_n\leq I_{\{z_n\not\in(b-B_c,b+B_c)\}},\quad n=1,2,\dots$$ is locally most powerful for testing $H_0:\theta=\theta_0$ vs $H_1:\theta>\theta_0$ if $b>0$ and $\phi_n=I_\{z_n\geq b+B_c \}$, $n=1,2,\dots$, and it is most powerful for testing $H_0:\theta=\theta_0$ vs $H_1:\theta<\theta_0$ if $b<0$ and  $\phi_n=I_\{z_n\leq b-B_c \}$, $n=1,2,\dots$. In both cases the optimality is in the class of all tests with the type I error probability and the average sample number not exceeding the corresponding values for $(\psi,\phi)$.

 In the case of $b=0$, both tests are locally most powerful, for the corresponding pair of hypotheses, in the class of all sequential tests whose average sample number does not exceed that of $(\psi,\phi)$.
 \end{remark}

\section{Appendix}
\subsection{Proof of Theorem \ref{t3.2}}
The proof is very close to the proof of Theorem 2.2 in \cite{novikov09b}.

First, the following lemma can be proved in exactly the same way as Lemma 5.1 in \cite{novikov09b}.
\begin{lemma} \label{l6.1}
Let, on a space with a $\sigma$-finite measure $\mu$, $F_1, F_2$ be some $\mu$-integrable functions
and $\phi$ some measurable function, such that
$$0\leq\phi(x)\leq 1. $$

Then
\begin{equation*}
\int(\phi(x)F_1(x)+(1-\phi(x))F_2(x))d\mu(x)\geq
\int\min\{F_1(x),F_2(x)\}d\mu(x)
\end{equation*}
with an equality if and only if
\begin{equation*}
I_{\{F_1(x)< F_2(x)\}}\leq\phi(x)\leq I_{\{F_1(x)\leq F_2(x)\}}
\end{equation*}
$\mu$-almost everywhere.
\end{lemma}

After this simple lemma, we can start with the proof of Theorem \ref{t3.2}.

Let us give to the left-hand
side of  (\ref{3.9}) the form
\begin{equation}\label{6.1}
b\alpha(\psi,\phi)-\dot\beta_0(\psi,\phi)=
\sum_{n=1}^\infty \int s_n^\psi
 \phi_n(bf_0^n-\dot f_0^n)
d\mu^n.
 \end{equation}

Applying Lemma \ref{l6.1} (with $F_2\equiv 0$) to each summand in
(\ref{6.1}) we immediately have:
 \begin{equation}\label{6.2}
 b\alpha(\psi,\phi)-\dot\beta_0(\psi,\phi)\geq \sum_{n=1}^\infty \int s_n^\psi
 \min\{0,bf_0^n-\dot f_0^n\}
d\mu^n.
 \end{equation}
 Let us note that the right-hand side of (\ref{6.2}) is finite: it follows from (\ref{6.1})
 by substituting  $\phi_n^\prime=I_{\{bf_0^n-\dot f_0^n<0\}}$ for $\phi_n$, $n=1,2,\dots$ in (\ref{6.1}).

Thus, there is an equality in (\ref{6.2}) if and only if each summand on the right-hand side of (\ref{6.1}) equals to the respective summand on the right-hand side of (\ref{6.2}).
And by Lemma \ref{l6.1} this happens if and only if
$\phi_n$ satisfies (\ref{3.10})
$\mu^n$-almost everywhere on $S_n^\psi$, for any $n=1,2,\dots$.
\subsection{Proof of Theorem \ref{t4.1}}
Using Lemma \ref{l6.1} instead of Lemma 5.1 of \cite{novikov09b} in the proof of Lemma 3.1 of \cite{novikov09b} we get
the following lemma, which takes over the major part of the proof of Theorem \ref{t4.1}.
\begin{lemma}\label{l6.2} Let $k$ be any integer non-negative  number, and let $$v_{k+1}=v_{k+1}(x_1,x_2,\dots,x_{k+1})$$ be
any $\mu^{k+1}$-integrable function. Then
$$
\sum_{n=1}^{k}\int s_n^\psi(cnf^n+l_n)d\mu^n
+\int t_{k+1}^\psi\left(c(k+1)f^{k+1}+v_{k+1}\right)d\mu^{k+1}
$$
\begin{equation}\label{6.3}
\geq\sum_{n=1}^{k-1}\int
s_n^\psi(cnf^n+l_n)d\mu^n+\int
t_{k}^\psi\left(ckf^{k}+v_{k}\right)d\mu^{k},
\end{equation}
with
\begin{equation*}
v_{k}=\min\{l_{k},cf^{k}+\int v_{k+1}d\mu(x_{k+1})\},
\end{equation*}
where, by definition,
$$\int v_{k+1}d\mu(x_{k+1})=\int v_{k+1}(x_1,\dots,x_{k+1})d\mu(x_{k+1}).$$
There is an equality in (\ref{6.3}) if and only if
\begin{equation*}
I_{\{l_{k}< cf^{k}+\int v_{k+1}d\mu(x_{k+1})\}}\leq\psi_{k}\leq I_{\{l_{k}\leq cf^{k}+\int v_{k+1}d\mu(x_{k+1})\}}
\end{equation*}
$\mu^{k}$-almost everywhere  on $T_{k}^\psi=\{(x_1,\dots,x_{k}): t_{k}^\psi(x_1,\dots,x_{k-1})>0\}$.
\end{lemma}

To start with the proof of Theorem \ref{t4.1}, let us first note  that, by definition, $Q_N^N(\psi)=L_N(\psi)$, and, by
Lemma \ref{l6.2},
\begin{equation*}
    L_N(\psi)\geq Q_{N-1}^N(\psi).
\end{equation*}

Also from Lemma \ref{l6.2} we easily get that
\begin{equation}\label{6.5}
    Q_{n+1}^N(\psi)\geq Q_{n}^N(\psi)
\end{equation}
for all $n=N-1,N-2,\dots,1$.

Thus, for any $k=1,\dots,N$
\begin{equation}\label{6.6}
   L_N(\psi)\geq Q_{k}^N(\psi).
\end{equation}
Obviously, there is an equality in (\ref{6.6}) if and only if there are equalities in all the inequalities in (\ref{6.5}), for all $n=k,k+1,\dots,N-1$. In turn, this happens, by the same Lemma \ref{l6.2}, if and only if
(\ref{4.5}) is satisfied $\mu^n$-almost everywhere on $T_n^\psi$ for all $n=k,k+1,\dots,N-1$.
\subsection{Proof of Theorem \ref{t4.2}}

The proof is very close to the proof of Theorem 3.2 in \cite{novikov09b}.
The same method is used in \cite{novikov09a} for multiple hypothesis testing.

Let $\psi\in\mathscr G$ be any stopping rule. By Lemma \ref{l4.2} for
any fixed $n\geq 1$
\begin{equation*}
    L(\psi)\geq Q_n(\psi).
\end{equation*}
In particular,
\begin{equation*}
    L(\psi)\geq Q_1(\psi)=c+R_0.
\end{equation*}

Passing the the limit in (\ref{6.5}), as $N\to\infty$, we have
\begin{equation}\label{6.9}
    Q_{n+1}(\psi)\geq Q_{n}(\psi),
\end{equation}
for any $n=1,2,\dots$, thus,
\begin{equation}\label{6.10}
    L(\psi)\geq Q_{n+1}(\psi)\geq Q_{n}(\psi)\geq c+R_0,
\end{equation}
for any $n=1,2,\dots$.
Supposing (\ref{4.16}), we have, by virtue of  Lemma \ref{l4.3}, that there are
equalities in all the inequalities in ({\ref{6.10}}). In particular, there is an equality in (\ref{6.9}), for any $n=1,2,\dots$.

Because, by the condition of Theorem \ref{t4.2}, $R_0>-\infty$, the integrals on both sides of (\ref{6.9}) are finite.
Applying Lemma \ref{l6.2}, we see that (\ref{4.17}) is fulfilled
$\mu^n$-almost everywhere on $T_n^\psi$, for any $n=1,2,\dots$.

(\ref{4.18}) now follows because
\begin{equation}\label{6.11}
Q_n(\psi)=L_n(\psi)+\int t_{n}^\psi(V_n-l_n)d\mu^n=c+R_0
\end{equation}
for any $n=1,2,\dots$, and $\lim_{n\to\infty}L_n(\psi)= L(\psi)=c+R_0$ by the conditions of the Theorem.

 The
``only if"-part of Theorem \ref{t4.2} is proved.

Let now $\psi$ satisfy (\ref{4.17}) $\mu^n$-almost everywhere on
$T_n^\psi$, for any $n=1,2,\dots$ and let (\ref{4.18}) hold for this $\psi$.

It follows from Lemma \ref{l6.2}  that
\begin{equation*}
    Q_n(\psi)=Q_{n-1}(\psi)=\dots=Q_1(\psi)=c+R_0
\end{equation*}
for any $n=1,2,\dots$.
It follows from (\ref{6.11}) and (\ref{4.18}) that $\lim_{n\to\infty}L_n(\psi)=c+R_0$. But $\psi\in \mathscr G$, so $\lim_{n\to\infty}L_n(\psi)=L(\psi)$, thus $L(\psi)=c+R_0=\inf_{\psi^\prime\in\mathscr G}L(\psi^\prime)$.
\subsection{Proof of Theorem \ref{t4.3}}

Let us first note that that if $\inf_{\psi\in\mathscr G} L(\psi;b,c)>-\infty$ for some $b\in \mathbb R$ and $c>0$, then $\inf_{\psi\in\mathscr G} L(\psi;b^\prime,c^\prime)>-\infty$ for all $b^\prime\in\mathbb R$ and $c^\prime\geq c$. Indeed,  if for any $k>-\infty$ there exists $(\psi,\phi)$ with $\psi\in\mathscr G$, such that $c^\prime E\tau_\psi+b^\prime\alpha(\psi,\phi)-\dot \beta_0(\psi,\phi)<k$, then
$$
c^\prime E\tau_\psi-\dot \beta_0(\psi,\phi)<k-\min\{0,b^\prime\},$$
and
$$
c E\tau_\psi+b\alpha(\psi,\phi)-\dot \beta_0(\psi,\phi)\leq c E\tau_\psi-\dot \beta_0(\psi,\phi)+\max\{0,b\}
$$
$$
<c^\prime E\tau_\psi-\dot \beta_0(\psi,\phi)+\max\{0,b\}\leq k-\min\{0,b^\prime\}+\max\{0,b\},
$$
so $\inf_{\psi\in\mathscr G} L(\psi;b,c)=-\infty$.

Thus, we can know that  $\inf_{\psi\in\mathscr G} L(\psi;b,c)>-\infty$, for all $b\in\mathbb R$, with just checking that
\begin{equation*}
h(c)=\inf_{\psi\in\mathscr G} L(\psi;0,c)>-\infty.
\end{equation*}

Let us note that $h:(0,\infty)\mapsto \mathbb R\cup \{-\infty\}$ is a concave function, as an infimum of a family of concave (linear) functions. In addition, it is obviously non-decreasing.
It easily follows from this, that either $h(c)>-\infty$ for all $c>0$ or $h(c)=-\infty$
for all $c>0$.

\subsection{Proof of Lemma \ref{l4.4}}

The first inclusion in (\ref{4.21}) is obvious.

Let $\psi\in \mathscr G_1$ be any stopping rule. Let us show that $L_N(\psi)\to L(\psi)$ as $N\to\infty$.
First, let us note that $L(\psi)$ is finite. This is because, on the one hand, by Lemma \ref{l4.3}, $L(\psi)>c+R_0>-\infty$. On the other hand, by (\ref{3.11}), $L(\psi)\leq E\tau_\psi<\infty$.

Now, by definition,
\begin{equation}\label{6.14} L(\psi)-L_N(\psi)=\sum_{n=N}^{\infty}\int
s_n^\psi(nf^n+l_n)d\mu^n
-\int t_N^\psi Nf^Nd\mu^N -\int t_N^\psi l_N d\mu^N \end{equation}

The first summand on the right-hand side of (\ref{6.14}) tends to $0$, as $N\to\infty$, because it is a tail of a converging series ($L(\psi)$).

The second summand on the right-hand side of (\ref{6.14}) tends to $0$ as well, because
$$ \int t_N^\psi N f^N d\mu^N = N E t_N^\psi=N P(\tau_\psi\geq N)\to 0,$$
as $N\to\infty$, since $E\tau_\psi<\infty$.

It remains to show that the third summand on the right-hand side of (\ref{6.14}) goes to $0$ as well.

To start with, let us note that
$$
\int t_N^\psi l_N d\mu^N=E_0\tau_\psi\left(\min\left\{0,b-{\dot f_0^N/ f_0^N}\right\}\right)
$$
(it is easy to see, using (\ref{2.1}), that $\dot f_0^n=0$ $\mu^n$-almost everywhere on $\{f_0^n=0\}$).
Therefore, using Schwarz' inequality we have
\begin{equation}\label{6.15} \left(\int t_N^\psi l_N d\mu^N\right)^2
   \leq
   E_0 (t_N^\psi)^2E_0\left(\min\left\{0,b-{\dot f_0^N/ f_0^N}\right\}\right)^2.\end{equation}
Because
$E_0(t_N^\psi)^2\leq E_0 t_N^\psi=P_0(\tau_\psi\geq N) $
and
$$
E_0(\min\left\{0,b-{\dot f_0^N/ f_0^N}\right\})^2\leq E_0(b-{\dot f_0^N/ f_0^N})^2$$
$$\leq 2 b^2+2E_0(\dot f_0^N/ f_0^N)^2\leq 2 b^2+2\gamma N
$$
for $N>N_0$ (by Assumption 3), we have from (\ref{6.15}) now that
$$
\left(\int t_N^\psi l_N d\mu^N\right)^2\leq 2P_0(\tau_\psi\geq N)(b^2+\gamma N)\to 0
$$
as $N\to\infty$, because $E_0\tau_\psi<\infty$.
\subsection{Proof of Theorem \ref{t5.2}}

Let $\psi\in \mathscr G$ be such that
\begin{equation}\label{6.16}\bar L(\psi;-b,c)=\inf_{\psi^\prime\in\mathscr G}\bar L(\psi^\prime;-b,c).\end{equation}

By Theorem \ref{t4.2}, it follows from (\ref{6.16}) that $\psi$ is defined by means of the functions $\bar l_n=\min\{0,-bf_0^n+\dot f_0^n\}$ and
$\bar V_n=\lim_{N\to\infty}\bar V_n^N$, $n=1,2,\dots$, where
$$
\bar V_n^N=\min\{\bar l_n,cf^n+\int \bar V_{n+1}^N d\mu(x_{n+1})\}
$$
$n=N-1,N-2,\dots,1$, being $V_N^N\equiv \bar l_N$.

Let us note that
\begin{equation}\label{6.17}
  \bar l_n= \min\{0,bf_0^n-\dot f_0^n\}-bf_0^n+\dot f_0^n=l_n-bf_0^n+\dot f_0^n.
\end{equation}
 Using this fact, it is easy to see that for any $N=1,2,\dots$
\begin{equation}\label{6.18}
\bar V_n^N=V_n^N-bf_0^n+\dot f_0^n,
\end{equation}
$\mu^n$-almost everywhere for any $n=N,N-1,\dots,1$.

Indeed, (\ref{6.18}) is satisfied for $n=N$ be virtue of (\ref{6.17}). Let us suppose now that (\ref{6.18}) is satisfied for some $n=k$. Then
$$
\bar V_{k-1}^N=\min\{\bar l_{k-1},cf^{k-1}+\int \bar V_k^N d\mu(x_k)\}$$
$$=\min\{\bar l_{k-1},cf^{k-1}+\int (V_k^N-bf_0^k+\dot f_0^k) d\mu(x_k)\}$$
$$=\min\{\bar l_{k-1},cf^{k-1}+\int V_k^Nd\mu(x_k)-bf_0^{k-1}+\dot f_0^{k-1} \},
$$ $\mu^{k-1}$-almost everywhere,
because $\int f_0^k(x_1,\dots,x_k)d\mu(x_k)= f_0^k(x_1,\dots,x_{k-1})$ and $\int \dot f_0^{k-1}(x_1,\dots,x_k)d\mu(x_k)= \dot f_0^{k-1}(x_1,\dots,x_{k-1})$ $\mu^{k-1}$-almost everywhere (the latter easily follows from (\ref{2.1})). Now, it follows from (\ref{6.17}) that (\ref{6.18}) is also satisfied for $n=k-1$, $\mu^{k-1}$-almost everywhere.

Now, passing to the limit, as $N\to\infty$, in (\ref{6.18}), we have
\begin{equation}\label{6.19}
    \bar V_n=V_n-bf_0^n+\dot f_0^n,
\end{equation}
$\mu^n$-almost everywhere.

Because, by Theorem \ref{t4.2},
$$
I_{\{\bar l_n< cf^n+\int\bar V_{n+1}d\mu(x_{n+1})\}}\leq\psi\leq I_{\{\bar l_n\leq cf^n+\int\bar V_{n+1}d\mu(x_{n+1})\}}
$$
$\mu^n$-almost everywhere on $T_n^\psi$, it follows from (\ref{6.17}) and (\ref{6.19}) that
\begin{equation}\label{6.20}
I_{\{l_n< cf^n+\int V_{n+1}d\mu(x_{n+1})\}}\leq\psi\leq I_{\{l_n\leq cf^n+\int V_{n+1}d\mu(x_{n+1})\}},
\end{equation}
$\mu^n$-almost everywhere on $T_n^\psi$, $n=1,2,\dots$.

In addition,
\begin{equation}\label{6.21}
\int t_n^\psi(V_n-l_n)d\mu^n=\int t_n^\psi(\bar V_n-\bar l_n)d\mu^n\to 0,
\end{equation}
as $n\to\infty$.

Let us show now that the problem of testing $H_0:\theta=\theta_0$ vs $H_1:\theta>\theta_0$ is finite.

It follows from (\ref{6.16}) that
\begin{equation}\label{56}
cE\tau_\psi-b \alpha(\psi,\phi)+\dot \beta_0(\psi,\phi)>k>-\infty
\end{equation}
for all $\psi\in \mathscr G$ and for all decision rules $\phi$. Let now $\bar \phi_n=1-\phi_n$, $n=1,2,\dots$.
Then $\alpha(\psi,\phi)=\sum_{n=1}^\infty E_0 s_n^\psi\phi_n=1-\sum_{n=1}^\infty E_0s_n^\psi\bar\phi_n$, and
$$
\dot \beta_0(\psi,\phi)=\sum_{n=1}^\infty \int \dot f_0^n s_n^\psi\phi_nd\mu^n=\sum_{n=1}^\infty  \int \dot f_0^ns_n^\psi d\mu^n-\sum_{n=1}^\infty \int \dot f_0^n s_n^\psi\bar \phi_nd\mu^n.
$$

It follows from (\ref{56}) now that
\begin{equation}\label{57}
cE\tau_\psi+b \alpha(\psi,\bar\phi)-\dot \beta_0(\psi,\bar\phi)+\sum_{n=1}^\infty\int \dot f_0^ns_n^\psi d\mu^n>k+b.
\end{equation}
Let $\tilde\phi_n\equiv 1$. Then, by Assumption 2,
$$
\left.\left(\sum_{n=1}^\infty\int f_\theta^n s_n^\psi\tilde\phi_nd\mu^n\right)_\theta^\prime\right|_{\theta=\theta_0}=
\sum_{n=1}^\infty\int \dot f_0^ns_n^\psi d\mu^n\leq 0,
$$
because $\sum_{n=1}^\infty\int f_\theta^n s_n^\psi\tilde\phi_nd\mu^n=P_\theta(\tau_\psi<\infty)=1$ for $\theta=\theta_0$, and is less than or equal to 1 for any other $\theta$.
We have from (\ref{57}) that
$$
cE\tau_\psi+b \alpha(\psi,\bar\phi)-\dot \beta_0(\psi,\bar\phi)>k+b,
$$
and therefore
\begin{equation}\label{6.24}
\inf_{\psi\in\mathscr G}L(\psi;b,c)\geq k+b >-\infty.
\end{equation}

It follows from (\ref{6.20}), (\ref{6.21}) and (\ref{6.24}), by Theorem \ref{t4.2}, that
\begin{equation}\label{6.25}
L(\psi;b,c)=\inf_{\psi^\prime}L(\psi^\prime;b,c).
\end{equation}

By analogy, it can be shown that  if $\psi$ satisfies (\ref{6.25}), then it satisfies (\ref{6.16}) as well.
\subsection{Proof of Theorem \ref{t5.4}}
The essential part of the proof is the use of Theorem \ref{t5.1} with $\mathscr G=\mathscr G_1$, where $\mathscr G_1=\{\psi:E_0\tau_\psi<\infty\}$ (see Lemma \ref{l4.4}, due to which the conditions of Lemma \ref{l4.3} are satisfied).

For the proof, we need some properties of the functions $\rho_c^n(z)$, $n=1,2,\dots$, and $\rho(z)$, $z\in\mathbb R$ (see (\ref{5.14}) and (\ref{5.16})).
\begin{lemma}\label{l6.3}Every $\rho_c^n:\mathbb R\mapsto \mathbb R$ is non-positive, concave, non-increasing, such that $\rho_c^n(z)+z$ is  non-decreasing with respect to $z$, and such that
\begin{equation*}
    \lim_{z\to\infty}(\rho_c^n(z)+z)=0\quad\mbox{and}\quad  \lim_{z\to-\infty}\rho_c^n(z)=0,
\end{equation*}
$n=1,2,\dots$.
\end{lemma}
{\em Proof.}
It is obvious that $\rho_c^0(z)=g(z)=\min\{0,-z\}$ has all the claimed in Lemma \ref{l6.3} properties.
If for some $k\geq 0$ $\rho_c^k$ has all these properties, let us show that $\rho_c^{k+1}$ does so as well.

By definition,
\begin{equation}\label{6.27}
    \rho_c^{k+1}(z)=\min\{g(z),c+\int \rho_c^k(z+r(x))f_{\theta_0}(x)d\mu(x)\}
\end{equation}
Obviously, $\rho_c^{k+1}$ is non-positive, concave and non-decreasing. Because of this, the integral on the right-hand side of (\ref{6.27}), by virtue of Lebesgue's monotone convergence theorem, goes to 0 as $z\to-\infty$. Thus, $\rho_c^{k+1}(z)\to 0$ as $z\to-\infty$.

We have further that
\begin{equation}\label{6.28}\rho_c^{k+1}(z)+z=\min\{\min\{z,0\},c+\int((z+r(x))+ \rho_c^k(z+r(x)))f_{\theta_0}(x)d\mu(x)\},\end{equation}
so it is non-decreasing with respect to $z$, and, by the monotone convergence theorem, goes to 0, as $z\to\infty$.
\rule{2mm}{2mm}\\
 It easily follows from Lemma \ref{l6.3} that $\rho_c(z)=\lim_{n\to\infty}\rho_c^n(z)$ possesses the same properties as $\rho_c^n$.
\begin{lemma}\label{l6.4}
The function $\rho_c:\mathbb R\mapsto \mathbb R$ is non-positive, concave, non-increasing, such that $\rho_c(z)+z$ is  non-decreasing with respect to $z$, and such that
\begin{equation}\label{6.29}
    \lim_{z\to\infty}(\rho_c(z)+z)=0\quad\mbox{and}\quad  \lim_{z\to-\infty}\rho_c(z)=0.
\end{equation}
\end{lemma}
{\em Proof.}
The only non-trivial thing to prove is (\ref{6.29}) -- other properties follow from the point-wise convergence.
To prove it, we start from
\begin{equation}\label{6.30}
    \rho_c(z)=\min\{g(z),c+\int\rho_c(z+r(x))f_{\theta_0}(x)d\mu(x)\},
\end{equation}
which follows from (\ref{6.27}) by the monotone convergence theorem, and
\begin{equation}\label{6.31}
    \rho_c(z)+z=\min\{\min\{z,0\},c+\int((z+r(x))+ \rho_c(z+r(x)))f_{\theta_0}(x)d\mu(x)\},
\end{equation}
which follows from (\ref{6.28}) in the same way.

Because $\rho_c(z)$ is non-increasing, $\lim_{z\to -\infty} \rho_c(z)=a_1$, so,
letting $z\to-\infty$ in (\ref{6.30}), we have that $a_1=\min\{0,c+a_1\}$, thus $a_1=0$.
Similarly, there exists $\lim_{z\to\infty}(\rho_c(z)+z)=a_2$. Passing to the limit, as $z\to\infty$, in (\ref{6.31}) we have $a_2=\min\{0,c+a_2\}$, thus $a_2=0$.
\rule{2mm}{2mm}

It easily follows from Lemma \ref{l6.4} that $g(z)-\rho_c(z)$ is a non-negative function tending to 0 as $z\to-\infty$ or $z\to \infty$, and such that
\begin{equation}\label{6.32}
    g(z)-\rho_c(z)\leq g(0)-\rho_c(0)=-\rho_c(0)\quad\mbox{for any}\quad z\in\mathbb R.
\end{equation}
Let $h_c(z)=\int \rho_c(z+r(x))f_{\theta_0}(x)d\mu(x)$ (see (\ref{6.30})). It follows from the Jensen inequality that $h_c(z)\leq \rho_c(z)$, $z\in\mathbb R$. In addition, obviously, $h_c(z)\to 0$ as $z\to -\infty$ and $h_c(z)+z\to 0$ as $z\to\infty$. Thus, $g(z)-h_c(z)$ is also a non-negative function tending to 0 as $z\to-\infty$ or $z\to\infty$ with a maximum reached at $z=0$.
In addition, it is easy to see that $g(z)-h_c(z)$ is convex as a function on $(-\infty,0]$, and that it is convex as a function on $[0,\infty)$.

Therefore, for any $0<c\leq -h_c(0)$ there are $A_c\leq 0$ and $B_c\geq 0$ such that
\begin{equation}\label{6.33}
    g(A_c)=c+h_c(A_c)\quad\mbox{and}\quad g(B_c)=c+h_c(B_c),
\end{equation}
and such that $g(z)<c+h_c(z)$ for $z<A_c$ or $z>B_c$, and $g(z)>c+h_c(z)$ for $z\in (A_c,B_c)$.
Because of this, (\ref{5.17}) is equivalent to (\ref{5.21}).
On the other hand, if $c> -h_c(0)$, then $c+h_c(z)>g(z)$ for all $z\in\mathbb R$. Thus, (\ref{5.17}) is equivalent to
$\psi_n=1$ in this case.

We have just proved that any $(b,c)$-generated stopping rule $\psi$ is as described immediately before Theorem \ref{t5.4}.

Let us show now that for any $(b,c)$-generated stoping rule $\psi$ it holds $E_0\tau_\psi<\infty$ (that is, $\psi\in\mathscr G_1$).
We have
$$
P_0(\tau_\psi\geq n)=E_0 t_n^\psi\leq P_{\theta_0}(t_n^\psi>0)
$$
$$
\leq P_{\theta_0}(z_k\in (b+A_c,b+B_c),\,\mbox{for any}\, k=1,2,\dots,n-1).
$$
Now, the finiteness of $E_0\tau_\psi$ follows by arguments of \cite{berk}, p. 376.

It is easy to see now that any $(b,c)$-generated stopping rule $\psi$ is regular, i.e. that  (\ref{5.19}) holds true. This is due to (\ref{6.32}),
because
$$
 0\leq\int t_n^\psi (g(z_n-b)-\rho_c(z_n-b))f_0^nd\mu^n\leq -\rho_c(0)\int t_n^\psi f_0^nd\mu^n
$$
$$
=-\rho_c(0) P_{\theta_0}(\tau_\psi\geq n)\to 0
$$
as $n\to\infty$.

Thus, all the conditions of Theorem \ref{t5.1} are satisfied, so the assertion of Theorem \ref{t5.4} follows.

\section{Acknowledgements}

The first author thanks  the National System of Investigators (SNI)
 of  CONACyT, Mexico for partial support for this work. The work of the first author is also partially supported by Mexico's CONACyT Grant no.
CB-2005-C01-49854-F.




\end{document}